\documentclass[amsmath,amssymb,prb,showkeys]{revtex4}

\newtheorem{theorem}{Theorem}[section]
\newtheorem{lemma}[theorem]{Lemma}
\newtheorem{proposition}{Proposition}[section]
\newtheorem{corollary}{Corollary}[section]

\numberwithin{equation}{section}

\newcommand{\blankbox}[2]

\usepackage{dcolumn}
\usepackage{bm}

\begin{document}

\title{Construct Weak Hopf Algebras By Using Borcherds Matrix}

\author{Wu Zhixiang }
\email{wzx@zju.edu.cn}
 \affiliation{Mathematics Department,\\
 Zhejiang University,Hangzhou, 310027, P.R.China}%


\begin{abstract}
We define a new kind quantized enveloping algebra of a generalized
Kac-Moody algebra ${\mathcal G}$ by adding a new generator $J$
satisfying $J^m=J$ for some integer $m$. We denote this algebra by
$wU_q^{\tau}({\mathcal G})$. This algebra is a weak Hopf algebra
if and only if $m=2,3$. In general, it is a bialgebra, and
contains a Hopf subalgebra. This Hopf subalgebra is isomorphic to
the usually quantum envelope algebra $U_q({\mathcal G})$ of a
generalized Kac-Moody algebra ${\mathcal G}$.
\end{abstract}

\keywords{generalized Kac-Moody algebras, irreducible
representation, weak Hopf algebra} \maketitle

\thanks{The author is sponsored by SRF for ROCS,SEM. He is also sponsored by the
Scientific Research Foundation of Zhejiang Provincial Education
Department.}

\maketitle

\section{Introduction}

 In his study of Monstrous moonshine$^{[3-5]}$, Borcherds introduced a new class of
 infinite dimensional Lie algebras called generalized Kac-Moody
 algebras. These generalized Kac-Moody algebra have a contravariant
 bilinear form which is almost positive definite. The fixed point algebra of any Kac-Moody algebra under a
 diagram automorphism is a generalized Kac-Moody algebra. A
 generalized Kac-Moody algebra can be regarded as a Kac-Moody algebra with imaginary simple
 roots. More explicitly, a generated Kac-Moody algebras is determined by a Borcherds-Cartan
 matrix $A=(a_{ij})_{(i,j)\in I\times I}$, where either $a_{ii}=2$, or $a_{ii}\leq 0$.
 If $a_{ii}\leq 0$, then the index $i$ is called imaginary, and the corresponding simple
 root $\alpha_i$ is called imaginary root. In this paper, the set $\{i\in
 I|a_{ii}=2\}$ is denoted by $I^+$. Set $I^{im}=I\setminus I^+$.
 The structure and the representation theory of generalized Kac-Moody
 algebras are very similar to those of Kac-Moody algebras, and
 many basic facts about Kac-Moody algebras can be extended to
 generalized Kac-Moody algebras. For example, the Kac-Weyl
 formula about an irreducible representation over a Kac-Moody algebra is generalized
 to a formula about an irreducible representation over a
 generalized Kac-Moody algebra as follows.
 $$chV(\lambda)=\frac{\sum_{w\in W}\sum_{F\subseteq
 T,F\perp \lambda}(-1)^{l(w)+|F|}e^{w(\lambda+\rho-s(F))}}{ \sum_{w\in
 W,F\subseteq T}(-1)^{l(w)+|F|}e^{w(\rho-s(F))}},$$
where $T$ is the set of all imaginary simple roots, $F$ runs all
over finite subsets of  $T$ such that any two elements in $F$ are
mutually perpendicular. We denote by $s(F)$ the sum of the roots
in $F$. We call the above formula Borcherds-Kac-Weyl formula.

On the other hand, many mathematicians are interested in
generalization of Hopf algebras, of which importance has been
recognized in both mathematics and physics. One way to do this is
to introduce a kind of weak coproduct such that $\Delta(1)\neq
1\otimes 1$ in Ref.{1}. The face algebras$^{[7]}$ and generalized
Kac algebras$^{[16]}$ are examples of this class of weak Hopf
algebras. Li and Duplij have defined and studied another kind of
weak Hopf algebras$^{[12]}$. A bialgebra
$(H,\mu,\eta,\Delta,\varepsilon)$ is called a weak Hopf algebra if
there is an anti-automorphism $T$ such that $T*id_H*T=id_H$ and
$id_H*T*id_H=T$, where $id_H$ is the identity map and $*$ is the
convolution product. Hopf algebras, and left of right Hopf$^{[
13-14]}$ algebras are weak Hopf algebras in this sense. In
presented paper a weak Hopf algebra is always mean the weak Hopf
algebra in this sense. The weak quantized enveloping algebras of
semi-simple Lie algebras are also weak Hopf algebras$^{[15]}$. Our
aim is to give more nontrivial examples of weak Hopf algebras.
Thanks to the definition of quantized enveloping algebra
$U_q({\mathcal G})$ associated a generalized Kac-Moody algebra
$\mathcal G$ defined in Ref. 7, we can also replace the group
$G(U_q({\mathcal G}))$ of grouplike elements by some regular
monoid as in Ref.14 and Ref. 15, we use a new generator $J$
instead of the projector in Ref.14 and Ref.15. Our new generator
$J$ satisfies $J^m=J$ for some integer $m\geq2$. By this way, we
obtain a subclass of bialgebra $wU_q({\mathcal G})$. This
bialgebra contains a sub-bialgebra $U^{\tau}(\mathcal G)$, which
is isomorphic to the a weak Hopf algebra in Ref. 15. Moreover, the
quotient algebra $wU_q({\mathcal
 G})/(1-J)$ is isomorphic to a sub-Hopf-algebra of the quantized enveloping algebra
 $U_q(\mathcal G)$ as Hopf algebras.
 As in the case of the classic quantum
group $U_q({\mathcal G})$, we try to determine irreducible
representation of $wU_q({\mathcal G})$.

Finally, let us outline the structure of this paper. In Section
II, we recall some basic facts related to the quantized enveloping
algebra of a generalized Kac-Moody algebra. In Section III, we
give the definition of $wU_q({\mathcal G})$. We study the
bialgebra structure of $wU_q({\mathcal G})$ in Section IV.  In the
final section, we study the irreducible representation of
$wU_q({\mathcal G})$.

\section{Notations and Preliminaries}
In this section, we fix notations and recall fundamental results
about generalized Kac-Moody algebras.

Let $I=\{1,\cdots,n\}$ or the set of positive integers, and
$A=(a_{ij})_{I\times I}$, a Borcherds-Cartan matrix, i.e., it
satisfies:

(1) $a_{ii}=2$ or $a_{ii}\leq 0$ for all $i\in I$,

(2) $a_{ij}\leq 0$ for all $i\not= j$,

(3) $a_{ij}\in {\bf Z}$,

(4) $a_{ij}=0$ if only if $a_{ji}=0$.

We say that an index $i$ is real if $a_{ii}=2$ and imaginary if
$a_{ii}\leq 0$. We denote $I^{+}=\{i\in I|a_{ii}=2\}$ and
$I^{im}=I-I^{+}$.

Kang considered the generalized Kac-Moody algebras associated with
Borcherds-Cartan matrices with charge$^{[9]}$ $${\bf m}=\{(m_i\in
{\bf Z}_{\geq 0})|i\in I, m_i=1\ for\ i\in I\}.$$ The charge $m_i$
is the multiplicity of the simple root corresponding to $i\in I$.
In this paper, we follow [11], and assume that $m_i=1$ for all
$i\in I$. However, we do not lose generality by this hypothesis.
Indeed, if we take Borcherds-Cartan matrices with some of the rows
and columns identical, then the generalized Kac-Moody algebras
with charge introduced in [9] can be recovered from the ones in
present paper by identifying the $h_is$ and $d_is$( and hence the
$\alpha_is$) corresponding to these identical rows and columns.

Moreover, we also assume that $A$ is symmetrizable; that is, there
is a diagonal matrix $D=diag\{s_i>0|i\in I\}$ such that $DA$ is a
symmetric matrix.

Let $P\ \check{}=(\oplus _{i\in I}{\bf Z}h_i)\oplus(\oplus_{i\in
 I}{\bf Z}d_i)$ be a free abelian group generated by the set
$\{h_i,d_i|i\in I\}$. This free abelian group is called the
co-weight lattice of $A$. The element $h_i$ in $\Pi
\check{}=\{h_i|i\in I\}$ is called a simple co-weight. We call
$\Pi \check{}$ the set of all simple co-weights. The space
${\mathcal H}={\bf Q}\otimes _{\bf Z}P\ \check{}$ over the
rational number field $\bf Q$ is said to be a Cartan subalgebra.
The weight lattice is defined to be $P:=\{\lambda\in {\mathcal
H}^*|\lambda(P\ \check{})\subseteq {\bf Z}\}$, where ${\mathcal
H}^*$ is the dual space of the Cartan subalgebra ${\mathcal
H}={\bf Q}\otimes _{\bf Z}P\ \check{}$. We denote  by $P^+$ the
set $\{\lambda\in P|\lambda(h_i)\geq 0, for\ every \ i\in I\}$ of
dominant integral weights.

Define $\alpha_i,\Lambda_i\in {\mathcal H}^*$ by
$$\alpha_i(h_j)=a_{ji},\  \  \   \alpha_i(d_j)=\delta_{ij}$$
$$\Lambda_i(h_j)=\delta_{ij}, \  \   \Lambda_i(d_j)=0.$$
Then $\alpha_i,i\in I$ are called simple roots of $A$. Let
$\Pi=\{\alpha_i|i\in I\}\subset P$ be the set of simple roots. The
free abelian group $Q=\oplus _{i\in I} {\bf Z}\alpha_i$ is called
the root lattice. Set $Q_+=\sum_{i\in I}{\bf Z}_{\geq 0}\alpha_i$
and $Q_-=-Q_+$. For any $\alpha\in Q_+$, we can write
$\alpha=\sum_{k=1}^n\alpha_{i_k}$ for $i_1,i_2,\cdots,i_n\in I$.
We set $ht(\alpha)=n$ and call it the height of $\alpha$.

Let $(.|.)$ be the bilinear form on $(\oplus_i({\bf
Q}\alpha_i\oplus {\bf Q}\Lambda_i))\times{\mathcal H}^*$ defined
by
$$(\alpha_i|\lambda)=s_i\lambda(h_i),(\Lambda_i|\lambda)=s_i\lambda(d_i).$$
Since it is symmetric on $(\oplus_i({\bf Q}\alpha_i\oplus {\bf
Q}\Lambda_i))\times(\oplus_i({\bf Q}\alpha_i\oplus {\bf
Q}\Lambda_i))$, one can extend this to a symmetric bilinear form
on ${\mathcal H}^*$. Then such a form is non-degenerated.

We always assume that ${\bf K}$ is a field of characteristic 0.
Let $q\in {\bf K}$ and $q_i=q^{d_i}$. It is assumed that $q_i\neq
\pm 1,0$ for all $i\in I$. For an indeterminant $\nu$ and an
integer $m$, let
$$[m]_{\nu}=\frac{\nu^m-\nu^{-m}}{\nu-\nu^{-1}},[m]!_{\nu}=[m]_{\nu}\cdots
[1]_{\nu},[0]_{\nu}=1,$$ and
$$\left[\begin{array}{c}
m\\
s
\end{array}\right]_{\nu}=\frac{[m]!_{\nu}}{[s]!_{\nu}[m-s]!_{\nu}}.$$
{\it {\bf Definition} The quantized enveloping generalized
Kac-Moody algebra $U_q({\mathcal G})$ associated with a
Borcherds-Cartan datum $(A,P$\v{},$P,\Pi$\v{}, $\Pi)$ is the
associated algebra with unit 1 over a field ${\bf K}$ of
characteristic 0, generated by the symbols $e_i$, $f_i$ $(i\in I)$
and $P\check{\quad}$ subject to the following defining relations:
$$q^0=1,q^hq^{h'}=q^{h+h'}\quad\forall h,h'\in P\check{},$$
$$q^he_iq^{-h}=q^{\alpha_i(h)}e_i,\qquad q^hf_iq^{-h}=q^{-\alpha_i(h)}f_i,$$
$$e_if_j-f_je_i=\delta_{ij}\frac {k_i-k_i^{-1}}{q_i-q_i^{-1}},\qquad where\quad k_i=q^{s_ih_i},$$
$$\sum_{r=0}^{1-a_{ij}}(-1)^r\left[\begin{array}{c}
1-a_{ij}\\
r
\end{array}\right]_i e_i^{1-a_{ij}-r}e_je_i^r=0\qquad if\quad a_{ii}=2,i\neq j,$$
$$\sum_{r=0}^{1-a_{ij}}(-1)^r\left[\begin{array}{c}
1-a_{ij}\\
r
\end{array}\right]_i f_i^{1-a_{ij}-r}f_jf_i^r=0\qquad if\quad a_{ii}=2,i\neq j,$$
$$ e_ie_j-e_je_i=f_if_j-f_jf_i=0,\qquad  if\quad a_{ij}=0.$$}

The quantum generalized Kac-Moody algebra $U_q({\mathcal G})$ has
a Hopf algebra structure with the comultiplication $\Delta$, the
counit $\varepsilon$, and antipode $S$ defined
by$$\Delta(q^h)=q^h\otimes q^h,$$ $$\Delta(e_i)=e_i\otimes
k_i^{-1}+1\otimes e_i,$$
$$\Delta(f_i)=k_i\otimes f_i+f_i\otimes
1,$$
$$\varepsilon(q^h)=1,\varepsilon(e_i)=\varepsilon(f_i)=0,$$
$$S(q^h)=q^{-h},S(e_i)=-e_ik_i,S(f_i)=-k_i^{-1}f_i$$
for all $h\in P\check{\quad}$ and $i\in I$.

Let $U_q^+({\mathcal G})$ and $U_q^-({\mathcal G})$ be the
subalgebras of $U_q({\mathcal G})$ generated by elements $e_i$ and
$f_i$ respectively, for $i\in I$, and let $U_q^0({\mathcal G})$ be
the subalgebra of $U_q({\mathcal G})$ generated by $q^h$($h\in
P\check{\quad}$). Then we have the triangular decomposition [1,8]

$$U_q({\mathcal G})=U_q^-({\mathcal G})\otimes U_q^0({\mathcal G})\otimes U_q^+({\mathcal G}).$$
Finally, let us use $U_q'(\mathcal G)$ to denote a subalgebra of
$U_q({\mathcal G})$ generated by $e_i,f_i,q^h$, where $h\in
\oplus_{i\in I} {\bf Z}s_ih_i\oplus\oplus_{i\in I} {\bf Z}d_i$. It
is obvious that $U_q'(\mathcal G)$ is a Hopf algebra.

\section{weak quantum algebras $wU_q^T({\mathcal G})$}
Let $m$ be a fixed positive integer. To generalize the
invertibility condition $k_ik_i^{-1}=1$ in $U_q({\mathcal G})$,
let us introduce some new generators $J$, $K_i$ and $\bar{K}_i$,
which subject the following relations:
\begin{eqnarray}J^{m-1}=K_i\bar {K}_i=\bar{K}_iK_i=D_i\bar {D}_i=\bar{D}_iD_i.
\end{eqnarray}Moreover, we assume that $J^m= J$ and
\begin{eqnarray}K_iJ^{m-1}=J^{m-1}K_i=K_i,\qquad \bar{K}_iJ^{m-1}=J^{m-1}\bar{K}_i=\bar{K}_i.
\end{eqnarray}\begin{eqnarray}
D_iJ^{m-1}=J^{m-1}D_i=D_i,\qquad
\bar{D}_iJ^{m-1}=J^{m-1}\bar{D}_i=\bar{D}_i.
\end{eqnarray}
 We call an element $E_i$ of type one if it satisfies
\begin{eqnarray}K_jE_i=q_i^{a_{ij}}E_iK_j,\quad \bar{K}_jE_i=q_i^{-a_{ij}}E_i\bar{K}_j .\end{eqnarray}
Similarly, if
\begin{eqnarray}K_jF_i=q_i^{-a_{ij}}F_iK_j,\quad \bar{K}_jF_i=q_i^{-a_{ij}}F_i\bar{K}_j,\end{eqnarray}
then $F_i$ is said to be type one. Suppose
\begin{eqnarray}K_jE_i\bar {K}_j=q^{a_{ij}}E_i.\end{eqnarray} Then we say
that $E_i$ is type zero. Similarly,  $F_i$ is type zero if it
satisfies the
following:\begin{eqnarray}K_jF_i\bar{K}_j=q^{-a_{ij}}F_i.\end{eqnarray}
\begin{proposition} $E_i$ ( resp. $F_i$ ) is type zero if and only if
$E_i$ is type one and $E_iJ^{m-1}=J^{m-1}E_i=E_i$ ( respectively,
$F_iJ^{m-1}=J^{m-1}F_i=F_i$).
\end{proposition}
{\it {\bf Proof}} \quad If $E_i$ is type zero, then  we obtain
from (III.6)
$$K_jE_i\bar{K}_jK_j=K_jE_iJ^{m-1}=q_i^{a_{ij}}E_iK_j.$$ On the other hand, since $J^{m-1}K_j=K_j$,
$$K_jE_i\bar{K}_j=K_jE_iJ^{m-1}\bar{K}_j=q_i^{a_{ij}}E_iK_j\bar{K}_j=q_i^{a_{ij}}E_iJ^{m-1}.$$
So $E_i=E_iJ^{m-1}$. Similarly, we can prove that
$E_i=J^{m-1}E_i$. Then
$$K_jE_i=K_jE_iJ^{m-1}=q_i^{a_{ij}}E_iK_j,$$and
$$E_i\bar{K}_j=J^{m-1}E_i\bar{K}_j=q_i^{a_{ij}}\bar{K}_jE_i.$$
That is, $E_i$ is type one.  On the other hand, if $E_i$ is type
one, and $E_iJ^{m-1}=J^{m-1}E_i=E_i$, then
$$q_i^{a_{ij}}E_i=q_i^{a_{ij}}E_iK_j\bar{K}_j=K_jE_i\bar{K}_j.$$
Similarly,we can prove the statement about $F_i$ is true. By now,
we complete the proof.

The types of $E_i$ and $F_i$ are denoted by $\kappa_i$,$\kappa'$
respectively. Let $\tau=(\{\kappa_i\}_{i\in I}|$
$\{\kappa_i')\}_{i\in I})$. By now, we can give the definition the
weak quantum algebra of type $\tau$ as follows:

{\it{\bf Definition} The type $\tau$ weak quantum algebra
$wU_q^{\tau}({\mathcal G})$ associated the generalized Kac-Moody
algebra ${\mathcal G}$ is an associated algebra with unit 1 over a
field ${\bf K}$ of characteristic 0, generated by  $J$,
$E_i,F_i(i\in I)$ and $K_i,D_i(i\in I)$ subjecting with the
following defining relations:
\begin{eqnarray}1\neq J^{m-1}=K_i\bar
{K}_i=D_i\bar{D}_i,\quad J^m=J,\end{eqnarray}
\begin{eqnarray}K_i\bar{K}_j=\bar{K}_jK_i,\quad K_iK_j=K_jK_i,\quad
\bar{K_i}\bar{K_j}=\bar{K}_j\bar{K}_i,\end{eqnarray}
\begin{eqnarray}D_i\bar{D}_j=\bar{D}_jD_i,\quad D_iD_j=D_jD_i,\quad
\bar{D_i}\bar{D_j}=\bar{D}_j\bar{D}_i,\end{eqnarray}
 \begin{eqnarray}D_i\bar{K}_j=\bar{K}_jD_i,\quad
K_iD_j=D_jK_i,\quad \bar{D_i}K_j=K_j\bar{D}_i, \quad
\bar{D_i}\bar{K}_j=\bar{K}_j\bar{D}_i,\end{eqnarray}\begin{eqnarray}D_iJ=JD_i,\quad
K_iJ=JK_i,\quad \bar K_iJ=J\bar K_i, \quad \bar{D_i}J =J\bar{D}_i,
\end{eqnarray}
\begin{eqnarray}D_iJ^{m-1}=D_i,\quad K_iJ^{m-1}=K_i,\quad J^{m-1}\bar{D_i}=\bar{D_i},\quad
J^{m-1}\bar{K}_i=\bar{K}_i,\end{eqnarray}
\begin{eqnarray} E_i\quad F_i\quad are \quad type\quad \tau,\end{eqnarray}

\begin{eqnarray}E_iF_j-F_jE_i=\delta_{ij}\frac {K_i-\bar{K}_i}{q_i-q_i^{-1}},\end{eqnarray}
\begin{eqnarray}\sum_{r=0}^{1-a_{ij}}(-1)^r\left[\begin{array}{c}
1-a_{ij}\\
r
\end{array}\right]_i E_i^{1-a_{ij}-r}E_jE_i^r=0\qquad if\quad a_{ii}=2,i\neq j,\end{eqnarray}
\begin{eqnarray}\sum_{r=0}^{1-a_{ij}}(-1)^r\left[\begin{array}{c}
1-a_{ij}\\
r
\end{array}\right]_i F_i^{1-a_{ij}-r}F_jF_i^r=0\qquad if\quad a_{ii}=2,i\neq j,\end{eqnarray}
\begin{eqnarray} E_iE_j-E_jE_i=F_iF_j-F_jF_i=0,\qquad  if\quad a_{ij}=0.\end{eqnarray}}

If $m=1$, and the Borcherds-Cartan matrix $A$ is symmetric, then
$wU_q^T({\mathcal G})=U_q({\mathcal G})$ provided that we identify
$K_i$ with $q^{h_i}$, $\bar {K}_i$ with $q^{-h_i}$, $D_i$ with
$q^{d_i}$ and $\bar{D}_i$ with $q^{-d_i}$. If $m=2$ and ${\mathcal
G}$ is a semisimple Lie algebra, then $wU_q^T({\mathcal G})$ has
been defined and studied by Yang in Ref.{[15]}. Notice that the
type zero was called type two by Yang.
\begin{lemma} $J^{m-1}$ is a center idempotent element of $wU_q^{\tau}({\mathcal G})$.
\end{lemma}
{\it {\bf Proof}}\quad If $E_i$ is type one, then
$J^{m-1}E_i=K_j\bar{K}_jE_i=E_iK_j\bar{K}_j=E_iJ^{m-1}$.
Similarly, we can prove that $F_iJ^{m-1}=J^{m-1}F_i$ provided
$F_i$ is type one. Hence this lemma follows from Proposition
III.1.

In the following corollary, the subalgebra of $U_q({\mathcal G})$
generated by $E_i,F_i$, $s_ih_i,d_i ( {i\in I} )$ is denoted by
$U_q'(\mathcal G)$. It is a Hopf subalgebra of $U_q(\mathcal G)$.
\begin{corollary} (1) $wU_q^{\tau}({\mathcal G})=wU_q^{\tau}({\mathcal G})J^{m-1}\oplus
wU_q^{\tau}({\mathcal G})(1-J^{m-1})$
is a direct sum of algebras.

(2) $wU_q^{\tau}({\mathcal G})/(1-J)$ is isomorphic to the algebra
$U_q'(\mathcal G)$.

(3) The subalgebra $wU_q^{\tau}({\mathcal G})_{m-1}$ of
$wU_q^{\tau}({\mathcal G})$, which is generated by
$E_iJ^{m-1},F_iJ^{m-1},K_i,\bar K_i,D_i,\bar D_i$ ($i\in I$), is
isomorphic to the algebra $U_q'(\mathcal G)$.

(4) Let $U^{\tau}({\mathcal G})$ be a subalgebra of
$wU_q^{\tau}({\mathcal G})$, which is generated by
$E_i,F_i,K_i,\bar K_i,D_i,\bar D_i$ ($i\in I$). Suppose ${\mathcal
G}$ is a semisimple Lie algebra. Then $U^{\tau}({\mathcal G})$ is
isomorphic to the algebra ${\mathcal M}^d_q(\mathcal G)$ in [15].
\end{corollary}
{\it {\bf Proof}}\quad The proof of (1) and (2) is easy. To prove
(3), let us define a map $\psi$ from $U_q'(\mathcal G)$ to
$wU_q^{\tau}({\mathcal G})J^{m-1}$ as follows.
$$\psi(q^{s_ih_i})=K_i,\quad \psi(q^{d_i})=D_i,\quad \psi(e_i)=E_iJ^{m-1},\quad \psi(f_i)=F_iJ^{m-1},\quad \psi(1)=
J^{m-1}.$$ Then one can show that $\psi$ is an algebra
homomorphism. Similarly we can define an algebra homomorphism a
map $\varphi$ from $wU_q^{\tau}({\mathcal G})J^{m-1}$ to
$U_q'(\mathcal G)$ as follows.
$$\varphi(K_i)=q^{s_ih_i},\quad \varphi(D_i)=q^{d_i},\quad \varphi(E_iJ^{m-1})=e_i,\quad \varphi(F_iJ^{m-1})=f_i,\quad
\varphi(J^{m-1})=1.$$ It is easy to check that $\varphi\psi=id$
and $\psi\varphi=id$. This proves (3).

To prove (4), let us use $J'$ to denote the $J$ in the algebra
${\mathcal M}^d_q(\mathcal G)$. Then we can define a mapping
$\theta$ from $U^{\tau}({\mathcal G})$ to ${\mathcal
M}^d_q(\mathcal G)$ as follows:
$$\theta(E_i)=E_i,\quad \theta(F_i)=F_i,\quad\theta(J^{m-1})=J'.$$ It is
easy to prove that $\theta$ is an isomorphism.

{\it {\bf Remark 1} By Proposition III.1, $wU_q^{\tau}({\mathcal
G})(1-J^{m-1})$ is generated by $1-J^{m-1}$, $E_i(1-J^{m-1})$ and
$F_i(1-J^{m-1})$, where $E_i$ and $F_i$ are type one. Since
(III.15) holds in $wU_q^{\tau}({\mathcal G})$,
$E_i(1-J^{m-1})F_j(1-J^{m-1})=F_j(1-J^{m-1})E_i(1-J^{m-1})$. If
all $\kappa_i=\bar{\kappa}_i=0$, then $wU_q^{\tau}({\mathcal
G})(1-J^{m-1})$ is isomorphic to  ${\bf K}$.}

\section{The bialgebra structure of $wU_q^{\tau}({\mathcal G})$}
The algebras $wU_q^{\tau}({\mathcal G})J^{m-1}$ and
$wU_q^{\tau}({\mathcal G})(1-J^{m-1})$ are denoted by $w$ and
$\bar {w}$ respectively in the following. By Corollary III.1,
$wU_q^{\tau}({\mathcal G})_{m-1}$ is isomorphic to the quantum
group $U_q(\mathcal G)$ provided $s_i=1$. Thus the
comultiplication and counit of $U_q(\mathcal G)$ transplant to the
algebra $wU_q^{\tau}({\mathcal G})_{m-1}$, and
$wU_q^{\tau}({\mathcal G})_{m-1}$ becomes a Hopf algebra.
Moreover, we can define three maps:
$$\Delta:wU_q^{\tau}({\mathcal G})\rightarrow wU_q^{\tau}({\mathcal G})\otimes wU_q^{\tau}({\mathcal G}),$$
$$\varepsilon:wU_q^{\tau}({\mathcal G})\rightarrow {\bf K},$$
$$T:wU_q^{\tau}({\mathcal G})\rightarrow wU_q^{\tau}({\mathcal
G}),$$as follows: \begin{equation}\Delta(K_i)=K_i\otimes K_i,\quad
\Delta(\bar{K}_i)=\bar{K}_i\otimes \bar{K}_i,
\end{equation}
\begin{equation}\Delta(D_i)=D_i\otimes D_i,\quad
\Delta(\bar{D}_i)=\bar{D}_i\otimes \bar{D}_i,
\end{equation}
\begin{equation}\Delta(J)=J\otimes J\end{equation}
\begin{eqnarray}\Delta(E_i)=\left\{\begin{array}{l}
1\otimes E_i+E_i\otimes K_i,\quad E_i\quad is \quad type\quad
one\\
J^{m-1}\otimes E_i+E_i\otimes K_i,\quad E_i\ is \ type\
zero,\end{array}\right.\end{eqnarray}
\begin{eqnarray}\Delta(F_i)=\left\{\begin{array}{l}
F_i\otimes 1+\bar{K}_i\otimes F_i,\quad F_i\quad is \quad
type\quad
one\\
F_i\otimes J^{m-1}+\bar{K}_i\otimes F_i,\quad F_i\ is \ type\
zero,\end{array}\right.\end{eqnarray}
\begin{eqnarray}
\varepsilon(K_i)=\varepsilon(\bar{K}_i)=1,\quad
\varepsilon(D_i)=\varepsilon(\bar{D}_i)=1,\quad \varepsilon(J)=1,
\end{eqnarray}
\begin{eqnarray}
\varepsilon(E_i)=\varepsilon(F_i)=0,\end{eqnarray}while the map
$T$ is defined as follows:
\begin{eqnarray}T(1)=1,\quad T(K_i)=\bar{K}_i,\quad
T(\bar{K}_i)=K_i ,\end{eqnarray}
\begin{eqnarray}T(J)=J,\quad T(D_i)=\bar{D}_i,\quad T(\bar{D}_i)=D_i ,\end{eqnarray}
\begin{eqnarray}T(E_i)=-E_i\bar{K}_i,\quad T(F_i)=-K_iF_i .\end{eqnarray}
Then we extend them to the whole $wU_q^{\tau}({\mathcal G})$. Thus
we obtain the following Lemma.

\begin{lemma}$wU_q^{\tau}({\mathcal G})$ is a bialgebra with
comultiplication $\Delta$ and counit $\varepsilon$.
\end{lemma}
{\it {\bf Proof}}\quad It can be shown by direct calculation that
the following relations
hold.$$\Delta(K_i)\Delta(\bar{K}_i)=\Delta(\bar{K}_i)\Delta(K_i)=\Delta(J^{m-1})=
\Delta(D_i)\Delta(\bar{D}_i)=\Delta(\bar{D}_i)\Delta(D_i) ,$$
$$\Delta(K_i)\Delta(\bar{K}_j)=\Delta(\bar{K}_j)\Delta(K_i),\quad
\Delta(D_i)\Delta(\bar{D}_j)=\Delta(\bar{D}_j)\Delta(D_i)
,$$
$$\Delta(K_i)\Delta({D}_j)=\Delta(D_j)\Delta({K}_i),\quad \Delta(\bar D_i)\Delta(\bar{K}_j)
=\Delta(\bar D_i)\Delta(\bar{K}_j),$$
$$\Delta(K_i)\Delta(\bar{D}_j)=\Delta(\bar D_j)\Delta({K}_i),\quad \Delta(D_i)\Delta(\bar{K}_j)
=\Delta(D_i)\Delta(\bar{K}_j),$$
$$\Delta(K_i)\Delta(J)=\Delta(J)\Delta({K}_i),\quad
\Delta(D_i)\Delta(J) =\Delta(D_i)\Delta(J),$$
$$\Delta(\bar K_i)\Delta(J)=\Delta(J)\Delta(\bar {K}_i),\quad
\Delta(\bar D_i)\Delta(J) =\Delta(\bar D_i)\Delta(J),\quad
\Delta(J^{m})=\Delta(J),$$
$$\Delta(J^{m-1}K_i)=\Delta({K}_i),\quad
\Delta(J^{m-1}\bar{K}_i)=\Delta(\bar{K}_i),$$
$$\Delta(J^{m-1}D_i)=\Delta({D}_i),\quad
\Delta(J^{m-1}\bar{D}_i)=\Delta(\bar{D}_i),$$
$$\varepsilon(K_i\bar{K}_j)=\varepsilon(K_i)\varepsilon(\bar{K}_j),\quad
\varepsilon(D_i\bar{D}_j)=\varepsilon(D_i)\varepsilon(\bar{D}_j),$$
$$\varepsilon(J^{m-1}K_i)=\varepsilon(K_i),\quad\varepsilon(J^{m-1}\bar{K}_j)=\varepsilon(\bar{K}_j),
$$
$$\varepsilon(K_i)\varepsilon(J)=\varepsilon(J)\varepsilon({K}_i),\quad
\varepsilon(D_i)\varepsilon(J) =\varepsilon(D_i)\varepsilon(J),$$
$$\varepsilon(\bar K_i)\varepsilon(J)=\varepsilon(J)\varepsilon(\bar {K}_i),\quad
\varepsilon(\bar D_i)\varepsilon(J) =\varepsilon(\bar
D_i)\varepsilon(J),\quad \varepsilon(J^{m})=\varepsilon(J),$$
$$
\varepsilon(J^{m-1}\bar{D}_i)=\varepsilon(D_i),\quad\varepsilon(J^{m-1}\bar{D}_j)=\varepsilon(\bar{D}_j),$$
$$
\varepsilon(D_iK_j)=\varepsilon(K_jD_i),\quad \varepsilon(D_i\bar
K_j)=\varepsilon(\bar K_jD_i),\quad \varepsilon(\bar
D_iK_j)=\varepsilon(K_j\bar D_i)$$
$$
\varepsilon(K_j)\varepsilon(E_i)=q_i^{a_{ij}}\varepsilon(E_i)\varepsilon(K_j),
\quad\varepsilon(F_i)\varepsilon(\bar{K}_j)=q_i^{a_{ij}}\varepsilon({F}_j)\varepsilon(\bar{K}_j),$$
$$
\varepsilon(E_i)\varepsilon(F_j)-\varepsilon(F_j)\varepsilon(E_i)=\delta_{ij}\frac
{\varepsilon(K_i)-\varepsilon(\bar{K}_i)}{q_i-q_i^{-1}}.$$ If
$E_i$ is type one, then $$\begin{array}{lll}
\Delta(K_j)\Delta(E_i)&=&(K_j\otimes K_j)((1\otimes E_i+E_i\otimes
K_i)\\
&=&K_j\otimes K_jE_i+K_jE_i\otimes
K_jK_i\\
&=&q_i^{a_{ij}}\Delta(E_i)\Delta(K_j).\end{array}$$ If $E_i$ is
type zero, then
$$\begin{array}{lll}
\Delta(K_j)\Delta(E_i)\Delta(\bar{K}_j)&=&(K_j\otimes
K_j)(J^{m-1}\otimes E_i+E_i\otimes K_i)(\bar{K}_j\otimes
\bar{K}_j)\\
&=&K_j\bar{K}_j\otimes K_jE_i\bar{K}_j+K_jE_i\bar{K}_j\otimes
K_jK_i\bar{K}_j\\
&=&q_i^{a_{ij}}\Delta(E_i).\end{array}$$
 Next we prove that
\begin{eqnarray}\Delta(E_i)\Delta(F_j)-\Delta(F_j)\Delta(E_i)=\delta_{ij}\frac{\Delta(K_i)-\Delta(\bar{K}_i)}{q_i-q_i^{-1}}
.\end{eqnarray} For any integers $0\leq r,s\leq m$, if
$J^{m-r}F_j=F_jJ^{m-r}$, $J^{m-r}\bar K_j=\bar K_jJ^{m-r}$,
$J^{m-s}K_i=K_iJ^{m-s}$ and $J^{m-s}E_i=E_iJ^{m-s}$, then
$$\begin{array}{l}
(J^{m-r}\otimes E_i+E_i\otimes K_i)(F_j\otimes
J^{m-s}+\bar{K}_j\otimes F_j)\\
\qquad -(F_j\otimes J^{m-s}+\bar{K}_j\otimes
F_j)(J^{m-r}\otimes E_i+E_i\otimes K_i)\\
=J^{m-r}F_j\otimes E_iJ^{m-s}+J^{m-r}\bar{K}_j\otimes
E_iF_j+E_iF_j\otimes K_iJ^{m-s}+E_i\bar{K}_j\otimes K_iF_j\\
-F_jJ^{m-r}\otimes J^{m-s}E_i-F_jE_i\otimes
J^{m-s}K_i-\bar{K}_jJ^{m-r}\otimes F_jE_i-\bar{K}_jE_i\otimes F_jK_i\\
=J^{m-r}\bar{K}_j\otimes (E_iF_j-F_jE_i)+ (E_iF_j-F_jE_i)\otimes
J^{m-s}K_i.
\end{array}
$$Thus (IV.11) holds for all $i,j$.

Finally, we prove that $\Delta$ satisfies the Quantum Serre
relations, i.e., the relations from (III.16) to (III.18).

From $a_{ij}=0$,  we obtain $E_iE_j=E_jE_i$, and
$K_iE_j=q_j^{a_{ji}}E_jK_i=E_jK_i$. Hence
$$\begin{array}{l} (J^{m-r}\otimes E_i+E_i\otimes
K_i)(J^{m-s}\otimes
E_j+E_j\otimes K_j)\\
\qquad -(J^{m-s}\otimes E_j +E_j\otimes
K_j)(J^{m-r}\otimes E_i+E_i\otimes K_i)\\
=J^{2m-r-s}\otimes E_iE_j+J^{m-r}E_j\otimes
E_iK_j+E_iJ^{m-s}\otimes K_iE_j+E_iE_j\otimes K_iK_j\\
-J^{2m-r-s}\otimes E_jE_i-J^{m-s}E_i\otimes E_j
K_i-E_jJ^{m-r}\otimes K_jE_i-E_jE_i\otimes K_jK_i\\
=J^{m-r}E_j\otimes E_iK_j+E_iJ^{m-s}\otimes
K_iE_j-J^{m-s}E_i\otimes E_j K_i-E_jJ^{m-r}\otimes K_jE_i\\
=0,
\end{array}$$
for  $ r=1,$ or $m-1,$ and $s=1$ or $m-1$. So
$$\Delta(E_i)\Delta(E_j)-\Delta(E_j)\Delta(E_i)=0.$$
Similarly, we can prove
$$\Delta(F_i)\Delta(F_j)-\Delta(F_j)\Delta(F_i)=0.$$
By now we have proven that $\Delta$ satisfies the relation (III.
18).

To prove that $\Delta$ satisfies the relation (III.16), we must
consider the following cases:

(1) Both $E_i$ and $E_j$ are type one.

(2) Only one of the $E_i$ and $E_j$ is type one.

(3) Both $E_i$ and $E_j$ are type zero.

For the case (3). Since $E_i$ is type zero,
$$\begin{array}{lll}
\Delta(E_i)&=&J^{m-1}\otimes E_i+E_i\otimes K_i\\
&=&(J^{m-1}\otimes 1)(1\otimes E_i+E_i\otimes K_i)\\
&=&(1\otimes E_i+E_i\otimes K_i)(J^{m-1}\otimes 1) \end{array}.$$
Let us introduce new notations. $(1\otimes E_i+E_i\otimes K_i)$ is
denoted by $\Delta'(E_i)$.  Set $s=1-a_{ij}$. Then
$$\begin{array}{l}
\sum_{r=0}^s(-1)^{r}\left[\begin{array}{c}s\\r\end{array}\right]\Delta(E_i)^{s-r}\Delta(E_j)\Delta(E_i)^r\\
=(J^{m-1}\otimes
1)^{(s+1)}\sum_{r=0}^s(-1)^{r}\left[\begin{array}{c}s\\r\end{array}\right]\Delta'(E_i)^{s-r}\Delta'(E_j)\Delta'(E_i)^r\\
=0, \end{array}$$ by the discussion in [8,pp67--68]. Hence
$\Delta$ satisfies (III.16) in this case. For the other cases, the
argument is more or less the same as case (3).

Similarly, we can prove that $\Delta$ satisfies(III.17). Therefore
$\Delta$ and $\varepsilon$ can be extended to an algebra morphism
from $wU_q^{\tau}({\mathcal G})$ to $wU_q^{\tau}({\mathcal
G})\otimes wU_q^{\tau}({\mathcal G})$, and from
$wU_q^{\tau}{\mathcal G})$ to ${\bf K}$, respectively.

It is easy to prove that \begin{eqnarray}(\Delta\otimes
1)\Delta(X)=(1\otimes \Delta)\Delta(X),\end{eqnarray}
\begin{eqnarray}(\varepsilon\otimes 1)\Delta(X)=(1\otimes\varepsilon
)\Delta(X)\end{eqnarray} for any
$X=E_i,F_i,K_i,\bar{K}_i,D_i,\bar{D}_i,J$. Since
$\Delta,\varepsilon$ are algebra morphisms, (IV.12) and (IV.13)
hold for any $X\in wU_q^{\tau} ({\mathcal G})$. By now we have
completed the proof.

Next we prove that the map $T$, defined by (IV.8),(IV.9),(IV.10),
is a weak antipode of the subbialgebra of the bialgebra
$wU_q^{\tau} ({\mathcal G})$ generated by $E_i,F_i,K_i,$
$\bar{K}_i,D_i,\bar{D}_i,J^{m-1}$. First we prove that $T$ can be
extended to an anti-automorphism of $wU_q^{\tau} ({\mathcal G})$.
It is easy to prove the following relations are
true.$$T(K_i)T(\bar {K}_j)=T(\bar {K}_j)T(K_i),\qquad T(D_i)T(\bar
{D}_j)=T(\bar {D}_j)T(D_i),$$
$$T(D_i)T(\bar {K}_j)=T(\bar {K}_j)T(D_i),\qquad T(K_i)T(\bar {D}_j)=T(\bar {D}_j)T(K_i),$$
$$T(\bar{D}_i)T(\bar {K}_j)=T(\bar {K}_j)T(\bar{D}_i),\qquad
T(J^{m-1})T(\bar{K}_i)=T(\bar{K}_i)\qquad
T(J^{m-1})T({K}_i)=T({K}_i),$$
$$T(K_iJ)=T(J)T({K}_i)=T(JK_i),\quad
T(D_iJ) =T(D_i)T(J) =T(JD_i),$$
$$T(\bar K_iJ)=T(J)T(\bar {K}_i)=T(J\bar K_i),\quad
T(\bar D_iJ) =T(\bar D_i)T(J)=T(J\bar D_i),\quad T(J^{m})=T(J),$$
$$T(J^{m-1})T(\bar{D}_i)=T(\bar{D}_i)\qquad T(J^{m-1})T({D}_i)=T({D}_i).$$
$$T(E_i)T(E_j)=T(E_j)T(E_i),\qquad T(F_i)T(F_j)=T(F_j)T(F_i),\quad
if \quad  a_{ij}=0.$$ If $E_i$ is type one, then
$$T(E_i)T(K_j)=-E_i\bar K_i\bar{K}_j=-q_i^{a_{ij}}\bar{K}_jE_i\bar K_i=q_i^{a_{ij}}T(K_j)T(E_i).$$
If $E_i$ is type zero, then
$$T(\bar{K}_j)T(E_i)T(K_j)=-K_jE_iK_i\bar{K}_j=-q_i^{a_{ij}}E_iK_i=q_i^{a_{ij}}T(E_i).$$
 Similarly, we can prove
 $$T(F_i)T(K_j)=q_i^{-a_{ij}}T(K_j)T(F_i)$$
 if $F_i$ is type one, and
 $$T(\bar{K}_j)T(F_i)T(K_j)=q_i^{-a_{ij}}T(F_i)$$ if $F_i$ is type
 zero. Moreover,
 $$\begin{array}{lll}
 T(F_j)T(E_i)-T(E_i)T(F_j)&=&K_j(F_jE_i)\bar{K}_i-E_i\bar{K}_iK_jF_j\\
 &=&q_j^{-a_{jj}}q_i^{a_{ij}}F_jE_iK_j\bar{K}_i-
 q_j^{-a_{jj}}q_j^{a_{ji}}F_jE_iK_j\bar{K}_i\\
 &=&\delta_{ij}\frac {\bar{K}_i-{K}_i}{q_i-q_i^{-1}}K_j\bar{K_i}\\
 &=&\delta_{ij}\frac {T(K_i)-T(\bar{K}_i)}{q_i-q_i^{-1}}
 \end{array}$$
 Similarly to [15, p.8], we can prove the following anti-relations to the
 quantum Serre relations hold.
$$\sum_{r=0}^s(-1)^{r}\left[\begin{array}{c}s\\r\end{array}\right]T(E_i)^{r}T(E_j)T(E_i)^{s-r}=0,\quad if\quad a_{ii}=2,$$
$$\sum_{r=0}^s(-1)^{r}\left[\begin{array}{c}s\\r\end{array}\right]T(F_i)^{r}T(F_j)T(F_i)^{s-r}=0\quad if\quad a_{ii}=2,$$
where $s=1-a_{ij}$.

From the above discussion, we get that $T$ is an anti-automorphism
of $wU_q^{\tau}(\mathcal G)$. Let $U^{\tau}$ be a subalgebra of
$wU_q^{\tau}(\mathcal G)$ generated by
$K_i,\bar{K}_i,D_i,\bar{D}_i,E_i,F_i,J^{m-1}$.
\begin{theorem} $T$ is a weak antipode of $U^{\tau}(\mathcal
G)$ and $U^{\tau}(\mathcal G)$ is a weak Hopf
algebra.\end{theorem}

{\it {\bf Proof}}\quad It is easy to verified the following
relations hold
$$(id*T*id)(X)=X,$$
$$(T*id*T)(X)=T(X),$$for
$X=K_i,\bar{K}_i,D_i,\bar{D}_i,E_i,F_i,J^{m-1}$.

Since $$id*T*id=(\mu\otimes1)\mu(id\otimes T\otimes
 id)(\Delta\otimes 1)\Delta,$$ $id*T*id$ is a linear automorphism of $wU_q^{\tau}(\mathcal
 G)$. To prove $(id*T*id)(X)=X,$ for any $X\in U^{\tau}(\mathcal
 G)$, we only need prove that \begin{eqnarray}(id*T*id)(xy)=xy\end{eqnarray} provided that
 $(id*T*id)(x)=x,$ and $y$ is one of the generators $K_i,\bar{K}_i,D_i,\bar{D}_i,E_i,$ $F_i,J^{m-1}$.
 Suppose $(\Delta\otimes 1)\Delta(x)=\sum x_{(1)}\otimes x_{(2)}\otimes x_{(3)}$.
Then $(\Delta\otimes 1)\Delta(xJ^{m-1})=\sum x_{(1)}J^{m-1}\otimes
x_{(2)}J^{m-1}\otimes x_{(3)}J^{m-1}$ and hence
$$id*T*id(xJ^{m-1})=\sum x_{(1)}T(x_{(2)})x_{(3)}J^{m-1}=xJ^{m-1}.$$ Similarly
$$\begin{array}{lll}
id*T*id(xE_i)&=&\sum
x_{(1)}J^{m-1}T(x_{(2)}J^{m-1})x_{(3)}E_i+\sum
x_{(1)}J^{m-1}T(x_{(2)}E_i)x_{(3)}K_i+\\
&&\sum
x_{(1)}E_iT(x_{(2)}K_i)x_{(3)}K_i\\
&=&\sum x_{(1)}T(x_{(2)})x_{(3)}E_i-\sum
x_{(1)}E_i\bar{K}_iT(x_{(2)})x_{(3)}K_i+\\
&&\sum
x_{(1)}E_i\bar{K}_iT(x_{(2)})x_{(3)}K_i\\
&=&xE_i ,\end{array}$$if $E_i$ is type zero. We can prove (IV.14)
is true for other generators of $U^{\tau}(\mathcal G)$. So
$id*T*id(x)=x$ for any $x\in U^{\tau}(\mathcal G)$ by induction.

Similarly, we can prove $T*id*T(x)=T(x)$ for any $x\in
U^{\tau}(\mathcal G)$. So $T$ is a weak antipode of
$U^{\tau}(\mathcal G)$, and $U^{\tau}(\mathcal G)$ is a weak Hopf
algebra.

\begin{corollary}  $wU_q^{\tau}({\mathcal G})$ is a weak Hopf algebra if and only if $m=2,3$.
Moreover, if $m=2,3$, then $wU_q^{\tau}({\mathcal G})$ is a
noncommutative and noncocommutative weak Hopf algebra with the
weak antipode $T$, but not a Hopf algebra.
\end{corollary}
{\it {\bf Proof}}\quad Since $m=2$, $J^2=J$. Then
$$id*T*id(J)=J^3=J=T(J)=T*id*T(J).$$ Similarly, if $m=3$, then
$J^3=J$. Thus $wU_q^{\tau}({\mathcal G})$
 is a weak Hopf algebra provided that $m=2,3$.

 If $wU_q^{\tau}({\mathcal G})$ is a weak Hopf algebra, then
$id*T*id(J)=J^3=J$. From this and $J^m=J$, we can obtain either
$J^2=J$ or $J^3=J$. Thus $m=2,3$ by our assumption.

 Suppose $wU_q^{\tau}({\mathcal G})$ is a Hopf algebra with antipode $S$. Then $S(J)J=1$. On the other
hand, since $J^m=J$, $(1-J^{m-1})J=0$ implies that $J^{m-1}=1$.
This is contradict to our assumption. So $wU_q^{\tau}({\mathcal
G})$ is not a Hopf algebra.

\begin{corollary} $U^{\tau}({\mathcal G})$ is a noncommutative
and noncocommutative weak Hopf algebra with the weak antipode $T$,
but not a Hopf algebra.  Moreover, $U^{\tau}({\mathcal G})J^{m-1}$
is isomorphic to $U_q'({\mathcal G})$ as Hopf algebras.
\end{corollary}
{\it {\bf Proof}}\quad It follows from Corollary III.1, Theorem
IV.2.
\begin{proposition}Suppose $m\geq 4$ and  $wU_q^{\tau}({\mathcal G})$ has a subalgebra $W$ containing $J^r$ for
some $1\leq r\leq m-1$. If $W$ is a weak Hopf algebra, then
\begin{eqnarray}r=\left\{\begin{array}{l}
m-1,\qquad if\quad  m \quad is\quad  even,\\
\frac 12(m-1),\qquad if\quad m\quad is\quad odd.\end{array}\right.
\end{eqnarray}
\end{proposition}
{\it {\bf Proof}}\quad Since $W$ is a weak Hopf algebra, $J^r=T*id
*T(J)=J^{3r}$. Hence $J^{3r}=J^r$ and $2r=k(m-1)$ for some integer
number $k$. If $m$ is even, $k=2$, then $r=m-1$. If $m$ is odd,
then $k=1,2$, and $r=m-1,\frac 12(m-1)$ respectively.

To simplify notations, let us use ${\mathcal B}_m$ to denote the
bialgebra $wU_q^{\tau}({\mathcal G})$, where $m$ is the minimal
integer number satisfying $J^m=J$.
\begin{proposition} Suppose $k=2,r=\frac12(m-1)$.
The subalgebra $W_3$ of $wU_q^{\tau}({\mathcal G})$ generated by
$E_i,F_i,K_i,\bar K_i,D_i,\bar D_i,J^r$ for $i\in I$ is isomorphic
to ${\mathcal B}_3$.

Suppose $k=1,r=m-1$. The subalgebra $W_2$ of
$wU_q^{\tau}({\mathcal G})$ generated by $E_i,F_i,K_i,\bar
K_i,D_i,\bar D_i,J^r$ for $i\in I$ is isomorphic to ${\mathcal
B}_2$.
\end{proposition}
{\it {\bf Proof}}\quad If $k=1$, and $r=m-1$, then $W_2$ is
isomorphic to ${\mathcal B}_2$ by Corollary IV.2.

If $m$ is odd, $k=2$, and $r=\frac 12(m-1)$, then we can prove
that $W_3$ is isomorphic to ${\mathcal B}_3$ similarly.

 Let $H$ be a coalgebra. The set of group-like elements of
$H$ is denoted by $G(H)$ in the next proposition.
\begin{proposition} $G(U^{\tau}({\mathcal G}))=G(U^{\tau}(\mathcal G)J^{m-1})\cup\{1\}.$
\end{proposition}
{\it {\bf Proof}}\quad If $g\in G(U^{\tau}({\mathcal G}))$, then
$g=gJ^{m-1}+g(1-J^{m-1})$. Let $g_1=gJ^{m-1},g_2=g(1-J^{m-1})$.
Then $g\otimes g=\Delta(g)=g_1\otimes g_1+g_1\otimes
g_2+g_2\otimes g_1 +g_2\otimes g_2$.  Since
$\Delta(g_1)=g_1\otimes g_1$ is a group-like element,
$\Delta(g_2)=g_1\otimes g_2+g_2\otimes g_1+g_2\otimes g_2$. So
$$\begin{array}{lll}(1\otimes \Delta)\Delta(g_2)&=&g_1\otimes
g_1\otimes g_2+g_1\otimes g_2\otimes g_1\\
& +&g_1\otimes g_2\otimes g_2+g_2\otimes g_1\otimes
g_1\\
&+&g_2\otimes g_1\otimes g_2+g_2\otimes g_2\otimes
g_1\\
&+&g_2\otimes g_2\otimes g_2.\end{array}$$ Then
\begin{eqnarray}\left\{\begin{array}{l}
(T*id*T)(g_2)=T(g_2)g_2T(g_2)=T(g_2),\\
(id*T*id)(g_2)=g_2T(g_2)g_2=g_2.\end{array}\right.
\end{eqnarray} Because $U^{\tau}(\mathcal G)(1-J^{m-1})$ is generated by
$E_i(1-J^{m-1}),F_j(1-J^{m-1}), 1-J^{m-1}$ and
$T(E_i(1-J^{m-1}))=T(F_j(1-J^{m-1}))=0$, $T(g_2)=k(1-J^{m-1})$ for
some $k\in {\bf K}$. From (IV.16), we obtain the following:
\begin{eqnarray}\left\{\begin{array}{l}
k^2g_2=k^2(1-J^{m-1})^2g_2=k(1-J^{m-1}),\\
kg_2^2(1-J^{m-1})=kg^2_2=g_2.\end{array}\right.
\end{eqnarray}If $k=0$, then $g=g_1\in G(U^{\tau}({\mathcal G}))$.
If $k\neq 0$, then $g_2=\frac 1k(1-J^{m-1})$. Thus
$$\frac 1k(1\otimes 1-J^{m-1}\otimes J^{m-1})=\frac1kg_1\otimes
(1-J^{m-1})+\frac1k(1-J^{m-1})\otimes
g_1+\frac1{k^2}(1-J^{m-1})\otimes (1-J^{m-1}).$$Multiplying by
$k(J^{m-1}\otimes 1)$ on the both sides of the above equation, we
get
$$J^{m-1}\otimes 1-J^{m-1}\otimes J^{m-1}=g_1\otimes (1-J^{m-1}).$$ Similarly, we
have $$1\otimes J^{m-1}-J^{m-1}\otimes J^{m-1}=(1-J^{m-1})\otimes
g_1.$$ Then
$$1\otimes 1-J^{m-1}\otimes J^{m-1}=J^{m-1}\otimes 1+1\otimes J^{m-1}-2J^{m-1}\otimes J^{m-1}+\frac1{k}(1-J^{m-1})\otimes
(1-J^{m-1}).$$Hence $$(1-J^{m-1})\otimes
(1-J^{m-1})=\frac1{k}(1-J^{m-1})\otimes (1-J^{m-1}).$$
Consequently, $k=1$. Notice that the set of group-like elements of
$U^{\tau}(\mathcal G)J^{m-1}$ is a monoid generated by
$K_i,\bar{K}_i,D_i,\bar{D}_i$ and $J^{m-1}$, and the elements from
this monoid are linearly independent over ${\bf K}$. So we get
$g_1=J^{m-1}$ from $1\otimes J^{m-1}-J^{m-1}\otimes
J^{m-1}=(1-J^{m-1})\otimes J^{m-1}=(1-J^{m-1})\otimes g_1$. Hence
$g=1$.

\begin{proposition} Suppose $\varphi$ is an automorphism of the
bi-algebra $wU^{\tau}_q(\mathcal G)$. Then $\varphi(J)=J^r$ for
some $1\leq r\leq m-1$ only if there exists $1\leq s\leq m-1$ such
that $rs \equiv 1mod(m-1)$. Moreover $\varphi(J^{m-1})=J^{m-1}$
and the restriction of $\varphi$ on $w$ ( resp. $\bar{w}$ )is an
isomorphism of $w$ (respect. $\bar{w}$).
\end{proposition}
{\it {\bf Proof}}\quad From $J^m=J$, we obtain
$\varphi(J)^m=\varphi(J)$. Thus
$\varphi(J)^mJ^{m-1}=\varphi(J)J^{m-1}$ is a group-like element in
$wU^{\tau}_q(\mathcal G)\subseteq U_q(\mathcal G)$. Suppose
$\varphi(J)=\varphi(J)J^{m-1}+x$, where
$x=\varphi(J)(1-J^{m-1})\in \bar{w}$.  Since
$T\varphi(J)=\varphi(T(J))=\varphi(J)$, $Tx=x$. Consequently,
$x=k(1-J^{m-1})$. Then
\begin{eqnarray}\begin{array}{ll}
\Delta(\varphi(J))&=(\varphi(J)J^{m-1}+x)\otimes(\varphi(J)J^{m-1}+x)\\
& =\varphi(J)J^{m-1}\otimes \varphi(J)J^{m-1}+k(1\otimes 1-J^{m-1}\otimes J^{m-1})\\
&=\varphi(J)J^{m-1}\otimes
\varphi(J)J^{m-1}+k\varphi(J)J^{m-1}\otimes (1-J^{m-1})\\
&+k(1-J^{m-1})\otimes\varphi(J)J^{m-1}
\end{array}
\end{eqnarray}From (IV.18), we obtain that
\begin{eqnarray}k(1\otimes 1-J^{m-1}\otimes J^{m-1})=k\varphi(J)J^{m-1}\otimes (1-J^{m-1})
&+k(1-J^{m-1})\otimes\varphi(J)J^{m-1}.\end{eqnarray} Hence
$$k(1-J^{m-1})=k\varphi(J)J^{m-1}(1-J^{m-1})+k(1-J^{m-1})\varphi(J)J^{m-1}=0.$$
So $k=0$ and $\varphi(J)=\varphi(J)J^{m-1}$. Notice that the
monoid of all group-like elements of $wU^{\tau}_q(\mathcal G)$ is
generated by $K_i,D_i,\bar K_i,\bar D_i$ and $J$. Since
$\varphi(J^m)=\varphi(J)^{m}=\varphi(J)$ $\varphi(J)=J^r$ for some
$1\leq r\leq m-1$.  Consequently $\varphi(J^{m-1})=J^{m-1}$, and
the restriction of $\varphi$ on $w$ is a Hopf algebra isomorphism
of $w$. Similarly, we can prove the restriction of $\varphi$ on
$\bar{w}$ is an algebra isomorphism of $\bar{w}$. Since $\varphi$
is an isomorphism, there exists $s$ such that $1\leq s\leq m-1$
and $$\varphi(J^s)=J^{rs}=J.$$Thus $rs\equiv1mod(m-1)$.

\noindent{\it{\bf Remark}\quad If  $r,s$ satisfy
$rs\equiv1mod(m-1)$, then it is easy to prove that the mapping
$\varphi$ defined as follows is an automorphism of the bialgebra
$wU^{\tau}_q(\mathcal G)$:
$\varphi(E_i)=E_i,\varphi(F_i)=F_i,\varphi(K_i)=K_i,\varphi(\bar
K_i)=\bar K_i,\varphi(D_i)=D_i,\varphi(\bar D_i)=\bar
D_i,\varphi(J)=J^r$. The inverse mapping $\psi$ of $\varphi$ is
given by the
following:$\psi(E_i)=E_i,\psi(F_i)=F_i,\psi(K_i)=K_i,\psi(\bar
K_i)=\bar K_i,\psi(D_i)=D_i,\psi(\bar D_i)=\bar D_i,\psi(J)=J^s.$

If the equation $rs\equiv1mod(m-1)$ has only trivial solution
$r=s=1$, for example, $m-1$ is a prime, then we can prove the
following:
\begin{corollary} Suppose $\mathcal G$ is a semi-simple Lie algebra, and the equation
$rs\equiv1mod(m-1)$ has trivial solution $r=s=1$. Then the
automorphism group of the bialgebra of $wU_q^{\tau}({\mathcal G})$
is the semi-direct product of $N$ and $H$, where $H$ is the group
of diagram automorphism, and $N$ is the group of diagonal
automorphism and it is a normal subgroup of the automorphism group
of $wU_q^{\tau}({\mathcal G})$.
\end{corollary}
{\it {\bf Proof}}\quad Similar to the proof of Theorem 5.1 in
Ref.15.

\noindent{\it{\bf Remark}\quad (1) If $\mathcal G$ is a semisimple
Lie algebra, then $U^{\tau}({\mathcal G})$ is isomorphic to the
weak quantum algebra defined by Yang in Ref. [15].

(2) Let ${\mathcal G}$ be a generated Kac-Moody algebra determined
by a Borcherds matrix $A$. Suppose $I^{im}\not=\phi$. Then
$\varphi_s$ for $1<s<m-1$ are automorphisms of
$wU_q^{\tau}({\mathcal G})$, where $\varphi_s$ defined as follows:
$\varphi_s(E_i)=E_i,\varphi_s(F_i)=F_i,\varphi_s(K_i)=K_i,\varphi_s(\bar
K_i)=$ $\bar K_i,\varphi_s(D_i)=J^{m-1-s}D_i,\varphi_s(\bar
D_i)=J^{s}\bar D_i$ for any $i\in I$.

\section{The representations of $wU_q^{\tau}({\mathcal G})$}

In this section, we try to determine the irreducible
representations of $wU_q^{\tau}({\mathcal G})$. Suppose $m-1$ is
invertible in the field $\bf K$ in this section.  Then $$e=\frac
1{m-1}\sum _{r=1}^{m-1}J^r$$ is a well-defined element. It is easy
to verify that $e^2=e$ and $J^{m-1}e=e$. Thus $e\in w$.

Suppose $V$ is a simple module over the bi-algebra
$wU_q^{\tau}({\mathcal G})$. Then $V=J^{m-1}V\oplus (1-J^{m-1})V$
is a direct sum of $wU_q^{\tau}({\mathcal G})$ modules. So either
$V=J^{m-1}V$, or $V=(1-J^{m-1})V$.

If $V=(1-J^{m-1})V$, then $Jv=0$ for any $v\in V$, and
$\bar{K}_iv=\bar{K}_iJ^mv=0$ for any $v\in V$.

If $V=J^{m-1}V$, then $J^{m-1}v=v$ for any $v\in V$. Suppose $v\in
V$ satisfying $K_iv=\lambda_iv$, then $\lambda_i\not=0$. In this
case, $K_i\bar{K}_iK_iv=\lambda_i^2\bar{K}_iv=\lambda_iv$. So
$\bar{K}_iv=\frac1{\lambda_i}v$.

By now we have completed the proof of the following proposition.
\begin{proposition} Let $V$ be simple $wU_q^{\tau}({\mathcal G})$-module.
Then either $J^{m-1}v=v$ for all $v\in V$, or $J^{m-1}v=0$ for any
$v\in V$. Suppose there exists an $i \in I$ such that
$K_{i}v=\lambda_{i}v$ for some nonzero vector $v$. Then
$\bar{K}_{i}v=\bar{\lambda}_{i}v$ for $\lambda_{i}$, where
$\bar{\lambda}_{i}=\left\{\begin{array}{ll} \lambda_{i}^{-1}& {\rm
if\ } \lambda_{i}\not=0;
\\0 &  {\rm if\ }\lambda_{i}=0.\end{array}\right.$ Moreover, $\lambda_i\neq 0$
if and only if $J^{m-1}v=v$.\end{proposition}

Let $V=J^{m-1}V$. Suppose $(F_ie-eF_i)V=(E_ie-eE_i)V=0$ for any
$i\in I$, then $eV$ is a module over $wU_q^{\tau}({\mathcal G})$.
Hence $V=eV\oplus (e-J^m)V$ is a sum of $wU_q^{\tau}({\mathcal
G})$ modules. So either $V=eV$,or $V=(e-J^{m-1})V$.

If $V=eV$, then $Jv=v$ for any $v\in V$. So $V$ can be viewed as a
module over $wU_q^{\tau}({\mathcal G})/(1-J).$ Notice that
$wU_q^{\tau}({\mathcal G})/(1-J)$ is isomorphic to $U'_q(\mathcal
G)$ by Corollary III.1(2). In this case, $V$ has been studied by
Kang $^{[9]}$. For example, the limit of highest weight simple
module is a highest weight simple module over the generalized
Kac-Moody algebra $\mathcal G$ with the same weight $\lambda$.
Then this simple module is unique determined by its formally
Borcherd-Kac-Weyl character formula ( see Section I).

If $V=(e-J^{m-1})V$, then $ev=0$ for any $v\in V$. Suppose a
nonzero element $v\in V$ satisfies $K_iv=\lambda v$ and $Jv=\gamma
v$. Then $\gamma+\cdots+\gamma^{m-1}=0$. So $\gamma^{m-1}=1$ and
$\gamma$ is a primitive $(m-1)$-th root of 1.

Suppose $V=(1-J^{m-1})V$ and $\bf K$ is an algebraically closed
field. Then $JV=0$ and $K_iV=\bar{K}_iV=0$ for any $i\in I$ by
Proposition V.1. Similarly, we can prove that $D_iV=\bar{D}_iV=0$
for all $i\in I$. Hence $E_iF_jV=F_jE_iV$ for all $i,j\in I$.
Moreover, $V$ can be viewed as a module over $\bar{w}$. Recall
that $\bar{w}$ is generated by $1-J^{m-1}$, $E_i(1-J^{m-1}),$ and
$F_j(1-J^{m-1})$, where $E_i,F_j$ are type one. Hence
$E_i(1-J^{m-1})F_j(1-J^{m-1})V=F_j(1-J^{m-1})E_i(1-J^{m-1})V$ for
all $i,j\in I$. In the following, we try to determine the
structure of $V$ in some special case.

Let $X_i=E_i(1-J^m),Y_i=F_i(1-J^m)$. Then every simple module $V$
over $\bar w$ is a module over the algebra generated by
$\{X_i,Y_j|$ for $i\in I_1$,$J\in I_2\}$, where $I_1=\{i\in I|E_i$
is type one$\}$, $I_2=\{j\in I|F_j$ is type one$\}$. The
generators $X_i,Y_j$ satisfy the following relation:
$$\begin{array}{c}
X_iY_j=Y_jX_i,\\
\sum_{r=0}^{1-a_{ij}}(-1)^r\left[\begin{array}{c}
1-a_{ij}\\
r
\end{array}\right]_i X_i^{1-a_{ij}-r}X_jX_i^r=0\qquad if\quad a_{ii}=2,i\neq
j,\\
\sum_{r=0}^{1-a_{ij}}(-1)^r\left[\begin{array}{c}
1-a_{ij}\\
r
\end{array}\right]_i Y_i^{1-a_{ij}-r}Y_jY_i^r=0\qquad if\quad a_{ii}=2,i\neq
j,\\
X_iX_j-X_jX_i=Y_iY_j-Y_jY_i=0,\qquad  if\quad a_{ij}=0.
\end{array}$$This simple module $V$ satisfies $JV=0$. From the above discussion, we obtain the following
result.
\begin{corollary}If $a_{ij}=0$ for any $i\in (I_1\cup I_2)\cap
I^+$, then every simple module over $\bar w$ is isomorphic to
$\bar w/M$, where $M$ is a maximal ideal of $\bar w$.
\end{corollary}
By Corollary V.1, the only simple over $\bar w$ is ${\bf
K}[x]/(p(x))$ if $|I_1\cup I_2|=1$, where $p(x)$ is an irreducible
polynomial in ${\bf K}[x]$. Suppose ${\bf K}$ is an algebraically
closed field. If $a_{ij}=0$ for any $i\in (I_1\cup I_2)\cap I^+$,
and $|I_1\cup I_2|=n$, then the simple module $V$ over $\bar w$ is
isomorphic to ${\bf K}[X_i,Y_j|i\in I_1,J_j\in
I_2]/(\{X_i-a_i,Y_j-b_j|i\in I_1,j\in I_2\})$ for some
$((a_i)_{i\in I_1},(b_j)_{j\in I_2})\in {\bf K}^n$.

\section*{ACKNOWLEDGMENT}

The author is sponsored by SRF for ROCS,SEM. He is also sponsored
by the Scientific Research Foundation of Zhejiang Provincial
Education Department. The author would like thank the referee for
his/her suggestion and help.


\begin{thebibliography}{99}
\bibitem{[1]} Aizawwa, N. and Isaac, P.S., Weak Hopf algebras corresponding to $U_q[sl_n]$,
J. Math. Phys. 44, 5250-5267(2003).

\bibitem{[2]} B\"ohm,G.,Nil,F.,and Szlach\'anyi,K. Weak Hopf algebras I,Integral theory and $C^*$-structure,
J.Algebra,221(1998),385-438.

\bibitem{[3]} R.E. Borcherds, Generalised Kac-Moody algebras, J. Algebra  115(1988), 501-512.

\bibitem{[4]} R.E. Borcherds, Monstrous moonshine and monstrous Lie superalgebras,
Invent. Math., 109(1992), 405-444.

\bibitem{[5]} R.E. Borcherds,Vertex algebras, Kac-Moody
algebras and the monster,Proc. Nat. Acad. Sci.
USA,83(1986),3068-3071.

\bibitem{[6]} Chin,W. and Musson,I.M., Corrigenda, the coradical filtration for quantized enveloping algebras,
J. Lond.Math. Soc. 61(2000),319-320.

\bibitem{[7]} Hayashi,T., An algebra related to the fusion rules of Wess-Zumino-Witten models,  Lett. Math. Phs. 22,
291-296(1991)

\bibitem{[8]} Jantzen, J.C., Lectures on Quantum Group (American Mathematical Society, Providence, RI, 1995) vol.6.

\bibitem{[9]}S.-J.Kang, Quantum deformations of generalised Kac-Moody
algebras and their modules,J. of Algebra 175(1995),1041-1066.

\bibitem{[10]} S.J.Kang, Quantum deformations of generalised Kac-Moody
algebras and their modules, J. of Algebra 175(1995), 1041-1066.

\bibitem{[11]} K. Jeong, S.-J.Kang, M. Kashiwara, Crystal bases for
quantum generalised Kac-Moody algebras and their modules, Proc.
London Math.Soc.(3)90(2205),395-438.

\bibitem{[12]} Li,F. and Duplij,S.,Weak Hopf algebras and singular
solutions of quantum Yang-Baxter equation, Commun. Math.
Phys.225(2002),191-217.

\bibitem{[13]} Nichols,W.D. and Taft,E.J. The left antipodes od a left Hopf algebra, Contemp.Math.
13(Ametrica Mathematical Society,Providence,RI,1982).

\bibitem{[14]} Wu Zhixiang, The Weak Hopf Algebras Related to Generalized
Kac-Moody Algebra, J.Math.Phs.(to appear).

\bibitem{[15]} Yang, S. ,Weak Hopf algebras corresponding to Cartan matrices,J.Math.Phys.46(2005),1-18.

\bibitem{[16]} Yamanouchi,T.,Duality for generalized Kacalgebras and a characterization of finite groupoid algebras,
J.algebra,163(1994),9-50.

\end{thebibliography}
\end{document}